\documentclass[12pt,a4paper]{article}
\usepackage{psfrag}
\usepackage{amssymb}
\usepackage[dvips]{graphicx}
\usepackage{epsfig}
\usepackage[usenames,dvipsnames]{color}

\newcommand{\stl}{\vspace{3mm}}

\newtheorem{defi}{\sc Definition}[section]
\newtheorem{theo}[defi]{\sc Theorem}
\newtheorem{prop}[defi]{\sc Proposition}

\newtheorem{lemm}[defi]{\sc Lemma}
\newtheorem{exem}[defi]{\sc Example}

\newenvironment{demo}[1][]{\noindent {\it Proof.} \protect\nopagebreak
       \rm #1}{\protect\nopagebreak $\square $ \par\stl}

\begin{document}

$$\mbox{\large  Logarithmic comparison theorem and ${\mathcal D}$-modules: an overview}$$

\stl

\begin{center}
 Tristan Torrelli\footnote{Laboratoire 
J.A. Dieudonn\'e, UMR
du CNRS 6621, Universit\'e de Nice Sophia-Antipolis,
Parc Valrose, 06108 Nice Cedex 2, France. {\it E-mail:} 
torrelli@math.unice.fr


 }
\end{center}

\bigskip

\ 

\stl

\noindent{\sc Abstract}. Let $D\subset X$ be a divisor
in a complex analytic manifold. A natural problem is to
determine when the de Rham complex of meromorphic 
forms on $X$ with poles along $D$ is quasi-isomorphic
to its subcomplex of logarithmic forms. In this  
mostly expository note, we recall the main results about 
this problem. In particular, we point out the relevance of the 
theory of $\cal D$-modules to this topic. 

\bigskip

\ 

\stl

\noindent{\bf\Large Introduction}

\stl

  Let $X$ be a complex analytic manifold of dimension $n\geq 2$.
Given a divisor $D\subset X$, we denote $j$ the natural inclusion
$X\backslash D\hookrightarrow X$. Let $\Omega^\bullet_X(\star
D)$ denote the complex of meromorphic forms on $X$ with poles
along $D$. From the Grothendieck Comparison Theorem \cite{gro}, the
de Rham morphism
$$\Omega^\bullet_X(\star D)\longrightarrow 
{\mathbf R}j_*{\mathbf C}_{X\backslash D}$$
is a quasi-isomorphism. In particular, if $X={\mathbf C}^n$, then
for each cohomology class $c\in H^p({\mathbf C}^n\backslash D,
{\mathbf C})$, there exists a differential form
$w\in\Omega^p_X(\star D)$ such that for any $p$-cycle $\sigma$ on
${\mathbf C}^n\backslash D$, one has $c(\sigma)=\int_\sigma w$.

\stl

  It is natural to ask what one can say about the form $w$. For
example, if $D$ is a complex submanifold then the order of the
pole of $w$ can be taken to be $1$. The question of the order of
the pole goes back to P.A. Griffiths \cite{Gri}. We recall that a
meromorphic form $w\in \Omega^p_X(\star D)$ is {\it logarithmic}
if $w$ and $dw$ have at most a simple pole along $D$; let
$\Omega^\bullet_X(\mbox{log}\,D)\subset \Omega^\bullet_X(\star D)$
denote the subcomplex of logarithmic forms with pole along $D$,  
introduced in full generality by K. Saito in \cite{saitlog}. In the initial case of
normal crossing divisors,  P. Deligne \cite{Del} proved that
the filtered morphism 
$(\Omega^\bullet_X(\mbox{log}\,D),\sigma)\hookrightarrow
(\Omega^\bullet_X(\star D),P)$ where $P$ is the
pole order filtration and $\sigma$ is induced by $P$,  
is a quasi-isomorphism compatible with
filtrations. This fact was crucial in order to defined a 
mixed Hodge structure on the cohomology of a
quasi-projective algebraic variety. 
Hence, one says that $D$ 
 satisfies the {\em logarithmic comparison theorem} if
 \begin{description}
   \item[\ \ LCT($D$)\,:] The inclusion
 $\Omega^\bullet_X(\mbox{log}\,D)\hookrightarrow \Omega^\bullet_X(\star D)$
  is a quasi-isomorphism.
\end{description}
A natural problem is therefore to find classes of divisors 
satisfying this condition, and also to understand its
meaning. Initiated by F.J. Castro-Jim\'enez, 
D. Mond and L. Narv\'aez-Macarro \cite{CMN}, 
this problem has been intensively studied these last years.
In this note, we gather together 
the main open questions\footnote{See the appendix.}
 and the main results. 
Essentially, they were obtained for hypersurfaces
with isolated singularities, hyperplane arrangements
and free divisors (see \S \ref{mainresul}). In this last case,
we recall the characterization in terms of 
$\cal D$-modules due to the Sevillian group around
F.J. Castro-Jim\'enez and L. Narv\'aez-Macarro 
(Theorem \ref{Dual}). Finally, we explain how enlightening 
this viewpoint is for the general study of the condition 
{\bf LCT($D$)} (see \S \ref{anal}).

\stl

\noindent{\bf Acknowledgements.} This research has 
been supported by a Marie Curie Fellowship of the European 
Community (programme
FP5, contract HPMD-CT-2001-00097). The author is
very grateful to the Departamento de \'Algebra, Geometr\'{\i}a y 
Topolog\'{\i}a (Universidad de Valladolid)
 for hospitality during the fellowship. Moreover, he is
happy to thank Francisco Jes\'us Castro-Jim\'enez, Antoine Douai,
  Luis Narv\'aez-Macarro and
Jos\'e Mar\'{\i}a Ucha-Enr\'{\i}quez  for judicious comments.

\section{Main results about {\bf LCT($D$)}}
\label{mainresul}

There are few families of divisors for which this condition 
{\bf LCT($D$)} has been studied. Indeed, it is difficult to work 
directly with the complex $\Omega^\bullet_X(\mbox{log}\,D)$ 
since we do not have in general a description of the logarithmic 
forms.

\subsection{The case of weighted homogeneous hypersurfaces with
an isolated singularity}

We recall that a polynomial
$h\in {\mathbf C}[x]={\mathbf C}[x_1,\ldots, x_n]$ is
{\it weighted homogeneous} of weight $d\in{\mathbf Q}^+$ for a system
$\alpha=(\alpha_1,\ldots,\alpha_n)\in({\mathbf Q}^{*+})^n$ if
$h$ is a (nontrivial) ${\mathbf C}$-linear combination of monomials
$x_1^{\gamma_1}\cdots x_n^{\gamma_n}$ with
$\sum_{i=1}^n \alpha_i\gamma_i=d$. In other words, we have the relation
$\chi(h)=dh$ where $\chi$ is the Euler-vector field
$\alpha_1x_1\partial_1+\cdots+\alpha_nx_n\partial_n$ associated
with $\alpha$.

\stl

 As usual, the case of weighted homogeneous
polynomials defining an isolated singularity at the origin provides
combinatorial formulas in terms of  the weights associated with the
 Jacobian algebra $A_h={\mathbf C}[x]/(h'_{x_1},\ldots,h'_{x_n})$. 
D. Mond and M. Holland \cite{HM} have obtained the 
following characterization:

\begin{theo}  \label{hachem}
   Assume that $n\geq 3$. Let $h\in{\mathbf C}[x]$
be a weighted homogeneous polynomial of degree $d$ for a system
$\alpha\in({\mathbf Q}^{*+})^n$.
Assume that $h$ defines an isolated singularity at the origin. Let
$D\subset {\mathbf C}^n$ be the hypersurface defined by $h$.
The following conditions are equivalent:
\begin{enumerate}
  \item The logarithmic comparison theorem holds for $D$.
  \item The link of $0$ in $D$ is a ${\mathbf Q}$-homology sphere.
  \item There is no weighted homogeneous element in  $A_h$
  whose weight belongs to the set
    $\{ k\times d-\sum_{i=1}^n\alpha_i\ ;\,1\leq k\leq n-2\}\subset
  {\mathbf Q}$.
\end{enumerate}
\end{theo}

In particular, the logarithmic comparison theorem does not hold in
general (see also Proposition \ref{ALCTWHIS}). For example, if 
$h=x_1^2+\cdots + x_n^2$ then we can take $d=2$, 
$\alpha_1=\cdots=\alpha_n=1$ and 
$A_h={\mathbf C}\overline{1}$.
Thus {\bf LCT($D$)} is satisfied if and only if $n=2$ or $n$ is odd.

\subsection{The case of hyperplane arrangements}

\label{arrang}

Let $D$ be a finite union of affine hyperplanes
$H$ in $X={\mathbf C}^n$, {\it i.e.} $H=\{\alpha_H=0\}$ where
$\alpha_H\in{\mathbf C}[x_1,\ldots,x_n]$ are polynomials of
degree one. We can associate with $D$
the ${\mathbf C}$-subalgebra of $\Omega^\bullet_X(\star D)$
generated by $1$ and the $1$-forms
$d\alpha_H/\alpha_H$. Let
$R^\bullet(D)$ denote this algebra of differential forms.
  It is well known that $R^\bullet(D)$ is isomorphic to the so-called
Orlik-Solomon algebra. Moreover

\begin{theo} {\em \cite{Brie}}
  For all $k\geq 0$, we have $R^k(D)\cong H^k(X\backslash D,{\mathbf C})$.
\end{theo}

On the other hand, we can consider the following complex of
${\mathbf C}$-vector spaces:
$0\rightarrow R^0(D)\stackrel{0}{\rightarrow} \cdots \stackrel{0}{\rightarrow}
  R^n(D)\rightarrow 0$ as a subcomplex of $\Omega^\bullet_X(\mbox{log}\,D)$.
Thus, a natural question is: does the logarithmic comparison theorem
hold for any hyperplane arrangement? This was conjectured by
H. Terao in \cite{Ter}. This is
true for tame arrangements (such as free arrangements, 
generic arrangements or  complex reflection arrangements) and when
$n\leq 4$ (see \cite{wy}). But in general, the question is still open.

\subsection{The case of free divisors}
\label{free}

Let ${\mathcal O}_X$ be the sheaf of holomorphic functions on
$X$. Given a divisor $D\subset X$, we will denote by
$h_D\in{\mathcal O}_{X,m}\cong{\mathcal O}= {\mathbf
C}\{x_1,\ldots, x_n\}$ a defining equation of $D$ at $m\in D$.

  A holomorphic vector field $v$ is {\em logarithmic} along
$D$ if for any point $m\in D$, $v(h_D)$ belongs to
$h_D{\mathcal O}_{X,m}$. Let $\mbox{Der}(-\mbox{log}\,D)$ denote
the (coherent) ${\cal O}_X$-module of logarithmic vector fields.
We recall a property studied by K. Saito in \cite{saitlog}.

\begin{defi}
  A divisor $D\subset X$ is {\em free at the point} $m\in D$ if
 $\mbox{\em Der}(-\mbox{\em log}\,D)_m$ is ${\cal O}_{X,m}$-free. 
It  is a {\em free divisor} if $\mbox{\em Der}(-\mbox{\em log}\,D)$ 
is  locally free.
\end{defi}

From the inclusions
$h_D \mbox{Der}({\mathcal O}_X)_m \subset  
\mbox{Der}(-\mbox{log}\,D)_m  \subset
 \mbox{Der}({\mathcal O}_X)_m$, the rank of
$\mbox{Der}(-\mbox{log}\,D)$ is also equal to $n$.

\begin{exem} {\em \label{exafree}
Free divisors appear in many different contexts.

$(i)$ Normal crossing divisors are free. Indeed, in
local coordinates such that $h_D=x_1\cdots x_p$, then
$$ \mbox{Der}(-\mbox{log}\,D)={\mathcal O}x_1\partial_1 \oplus
\cdots \oplus {\mathcal O}x_p\partial_p \oplus
{\mathcal O}\partial_{p+1} \oplus
\cdots \oplus {\mathcal O}\partial_n.$$

$(ii)$ Plane curves are free (K. Saito \cite{saitlog}).

$(iii)$ Complex reflection arrangements are free (H. Terao \cite{Tera}).
For example, the braid arrangement,  defined by
$\prod_{1\leq i<j\leq n}(x_i-x_j)$ in ${\mathbf C}^n$,  is free.

$(iv)$ The discriminant of a versal deformation of
an isolated complete intersection singularity is a free divisor 
(see \cite{saito}; \cite{Loo}, Corollary 6.13; \cite{Alek}).}
\end{exem}

Most of the known results about the logarithmic comparison theorem
 have been  obtained for free divisors. Indeed, the logarithmic de Rham 
complex $\Omega^\bullet_X(\mbox{log}\,D)$
is also explicit. More precisely, because of the duality between
 $\Omega^1_X(\mbox{log}\,D)$ and
$\mbox{Der}(-\mbox{log}\,D)$, $\Omega^1_X(\mbox{log}\,D)$
is also a free ${\mathcal O}_X$-module and we have 
$\Omega^q_X(\mbox{log}\,D)=
\bigwedge^q \Omega^1_X(\mbox{log}\,D)$ for $1\leq q\leq n$ (see 
\cite{saitlog}). 
Firstly, we have the following characterization for plane curves:

\begin{theo}{\em \cite{CCMN}} \label{CCMN}
  If $D\subset X={\mathbf C}^2$ is a plane curve, then
 the logarithmic comparison theorem holds if and only if
$D$ is locally weighted homogeneous.
\end{theo}

This last condition means that for all $m\in D$, there exists an
analytic change of coordinates $\phi$ such that $h_D\circ \phi$
is a weighted homogeneous polynomial; for example, weighted 
homogeneous hypersurfaces with an isolated singularity and 
hyperplane arrangements are locally weighted homogeneous. 
This unusual condition is the suitable one in this context for doing 
inductions on the dimension of $D$ (see the proof of Proposition 
\ref{BetLCT} for example). More generally, we have

\begin{theo} {\em \cite{CMN}} \label{CMN}
   Let $D\subset X$ be a locally weighted homogeneous free divisor.
Then the logarithmic comparison theorem holds for $D$.
\end{theo}

\noindent Among the free divisors in Example \ref{exafree}, the one 
given in $(i)$, $(iii)$ and some\footnote{For more details, 
see \cite{CMN}.} of  $(iv)$ are locally weighted homogeneous. 

\stl

The converse is false in general. For example, 
$h=x_1x_2(x_1+x_2)(x_1+x_2x_3)$ defines a free divisor 
$D\subset {\bf C}^3$ such that
  {\bf LCT($D$)} is true and $h$ is not weighted homogeneous (see 
\cite{CCMN}, \S 4). Meanwhile, $h$ belongs to the ideal of its partial 
derivatives. In other words, there exists locally a vector field $v$ such 
that $v(h_D)=h_D$;
  one says sometimes that $h$ is  {\it Euler-homogeneous}. 
In fact, we have no example
of a free divisor $D=V(h)$ verifying {\bf LCT($D$)} which is not
Euler-homogeneous. This is true for a Koszul-free divisor (see Definition
\ref{DKF}, Theorem \ref{carkosfree}); moreover, M. Granger and M. Schulze
\cite{schul} have obtained the following result:

\begin{theo}
 Let $D=V(h)\subset X={\bf C}^3$ be a free
divisor. If the logarithmic comparison theorem
holds for $D$, then $h$ is Euler-homogeneous.
\end{theo}

For $n\geq 4$, this question is still open (see \cite{CCMN}, Conjecture 1.4)
and it can be extended for a general divisor.

\section{A differential viewpoint for free divisors}

Here we recall how the condition {\bf LCT($D$)} may be interpreted 
in terms of ${\cal D}_X$-modules for free divisors $D\subset X$, 
as it was initiated by F.J. Calder\'on-Moreno in \cite{3}.

\subsection{Preliminaries}

Given a complex analytic manifold $X$ of dimension $n\geq 2$, we
denote $\Omega^\bullet_X$ the complex of holomorphic differential
forms on $X$ and $({\cal D}_X, F_\bullet )$ the sheaf of linear differential
operators with holomorphic coefficients filtered by order. 
Locally at a point $m\in X$, we have
${\cal O}_{X,m}\cong{\mathcal O}={\mathbf C}\{x_1,\ldots,x_n\}$ and 
${\cal D}_{X,m}\cong{\mathcal D}={\mathcal O}\langle
\partial_1,\ldots,
\partial_n\rangle$; moreover we identify $\mbox{gr}^F{\mathcal D}$
with ${\mathcal O}[\xi]={\cal O}[\xi_1,\ldots,\xi_n]$.

\stl

The  so-called Riemann-Hilbert correspondence
of Z. Mebkhout and M. Kashiwara (see \cite{K4}, \cite{M2}, \cite{ZZ}) 
asserts that there is an equivalence of categories between the
category $hr({\mathcal D}_X)$ of (left) regular holonomic
${\cal D}_X$-modules and the one of perverse sheaves
$Perv_X({\mathbf C})$ on $X$ {\em via} the de Rham functor
\begin{eqnarray*}
  hr({\mathcal D}_X) & \longrightarrow & Perv_X({\mathbf C}) \\
  {\mathcal M} & \longmapsto & \mbox{DR}({\mathcal M})=\Omega^\bullet_X
\otimes_{{\mathcal O}_X}{\mathcal M}.
\end{eqnarray*}
Roughly speaking, a perverse sheaf  on $X$ is a special type of
complex of sheaves on $X$ whose cohomology groups are
constructible in ${\mathbf C}$-vector spaces of finite dimension
on a stratification of $X$. For example, ${\mathcal O}_X$ is
regular holonomic and $\mbox{DR}({\mathcal
O}_X)=\Omega^\bullet_X$ is quasi-isomorphic to the constant sheaf
${\mathbf C}_X$  by the Poincar\'e lemma.

\subsection{On the perversity of $\Omega^\bullet_X(\mbox{log}\,D)$}

Given a divisor $D\subset X$, we consider the sheaf
${\mathcal O}_X(\star D)$ of meromorphic functions with
poles along $D$. As ${\mathcal O}_X(\star D)$ is
regular holonomic, the meromorphic de Rham complex
$\mbox{DR}({\mathcal O}_X(\star D))={\Omega^\bullet_X(\star D)}$
is a perverse sheaf too. Thus it is natural to investigate conditions
on $D$ in order to get the perversity of
$\Omega^\bullet_X(\mbox{log}\, D)$. In the case of free
divisors, this question was studied by F.J. Calder\'on-Moreno
and L. Narv\'aez-Macarro in \cite{3}, \cite{Dual}. They obtained
the following characterization:

\begin{theo} \label{Perv}
   Let $D\subset X$ be a free divisor. Then the logarithmic
complex $\Omega^\bullet_X(\mbox{\em log}\,D)$ is
perverse if and only if the following conditions
are satisfied:
\begin{enumerate}
  \item the complex ${\mathcal D}_X
\otimes^L_{{\mathcal V}^D_0({\mathcal D}_X)}{\mathcal O}_X$
  is concentrated in degree $0$;
  \item the ${\cal D}_X$-module 
${\mathcal D}_X\otimes_{{\mathcal V}^D_0({\mathcal D}_X)}{\mathcal O}_X$
is holonomic.
\end{enumerate}
\end{theo}

\noindent From \cite{7}, we say also that $D$ is of {\em Spencer type}. 
Here ${\mathcal V}^D_0({\mathcal D}_X)\subset{\mathcal D}_X$ 
is the coherent sheaf of rings\footnote{For some results about this sheaf 
for a general divisor, see \cite{Schu}.} 
of logarithmic operators \cite{3} (that is, $P\in{\mathcal D}_X$
such that locally  $P\cdot (h_D^k)\subset h_D^k{\mathcal O}$ for any
integer $k$). 
Let us notice that this condition $1$ has no clear meaning. Thus, the
 problem is now to find (geometrical) criteria on a free divisor to
be of Spencer type (see \cite{Dual}, \S 5). The only known condition
is to be a Koszul-free divisor (see \cite{3}, Theorem 4.2.1).

\begin{defi} \label{DKF}
   A free divisor $D\subset X$ is {\em Koszul-free} if there
exists locally a  basis\footnote{In fact, every basis of 
$\mbox{Der}(-\mbox{log}\,D)$
satisfies this property when $D$ is Koszul-free.}  $\{\delta_1,\ldots,\delta_n\}$ of
$\mbox{\em Der}(-\mbox{\em log}\,D)$ such that the sequence of
principal symbols $(\sigma(\delta_1),\ldots, \sigma(\delta_n))$ is
$\mbox{\em gr}^F{\cal D}$-regular.
\end{defi}

\noindent  For example, plane curves and locally weighted homogeneous
free divisors are Koszul-free \cite{saitlog}, \cite{CN}, \cite{LQH}. This notion
means also that the free divisor $D$ has a holonomic stratification
in the sense of K. Saito (see \cite{saitlog}, \S 3; \cite{Bru}, Proposition 6.3; 
\cite{3}, Corollary 1.9). 

\stl

Finally, we make the remark that being Koszul-free is not necessary for
the perversity of $\Omega^\bullet_X(\mbox{log}\,D)$; for example, 
the free divisor $D\subset{\bf C}^3 $ defined by 
$x_1x_2(x_1+x_2)(x_1+x_2x_3)$ is not Koszul-free but 
$\Omega^\bullet_X(\mbox{log}\,D)$ 
is perverse \cite{3}, \cite{CCMN}. On the other hand, the complex 
$\Omega^\bullet_X(\mbox{log}\,D)$
is not perverse for any free divisor. For example, the divisor
$(x_1^5+x_2^4+x_1^4x_2)(x_1+x_2x_3)=0$  
is free but  not of Spencer type (see \cite{Dual}, \S 5).

\subsection{A differential characterization of {\bf LCT($D$)}}

Let us now give a differential analogue of condition {\bf LCT($D$)}.
 F.J. Castro-Jim\'enez and J.M. Ucha-Enr\'{\i}quez
began work on this problem in \cite{CUE} \& \cite{7}, and they 
obtained a characterization for free divisors of Spencer type 
\cite{CUexp}. For a general free divisor, we have the following 
generalization\footnote{In fact, they obtained this result 
for any integral logarithmic connection (see \cite{Dual}, Theorem 4.1).} 
due to F.J. Calder\'on-Moreno and L. Narv\'aez-Macarro \cite{Dual}:
\begin{theo}  \label{Dual}
   Let $D\subset X$ be a free divisor. Then the
inclusion
 $$\Omega^\bullet_X(\mbox{\em log}\,D)
\hookrightarrow \Omega^\bullet_X(\star D)$$
is a quasi-isomorphism if and only if the following conditions
are satisfied:
\begin{enumerate}
  \item the complex ${\mathcal D}_X
\otimes^L_{{\mathcal V}^D_0({\mathcal D}_X)}{\mathcal O}_X(D)$
  is concentrated in degree $0$;
  \item the natural morphism
$$\varphi_D:{\mathcal D}_X
\otimes_{{\mathcal V}^D_0({\mathcal D}_X)}{\mathcal O}_X(D) \longrightarrow
{\mathcal O}_X(\star D)$$
is an isomorphism.
\end{enumerate}
\end{theo}
Here ${\mathcal O}_X(D)$ denotes the ${\mathcal
V}^D_0({\mathcal D}_X)$-module of meromorphic functions with
 at most a simple pole along $D$. Unfortunately, this characterization 
is no more explicit than condition  {\bf LCT($D$)}; meanwhile, 
condition $1$ is verified by Koszul-free divisors.

\stl

A key point in
the proofs of Theorems \ref{Dual} and \ref{Perv} is the
relationship between the duals of any integral logarithmic
connection over the base ring ${\cal D}_X$ and 
${\cal V}^D_0({\cal D}_X)$ (see \cite{7}; \cite{Dual}, \S 3). 
In particular, 
the holonomic ${\cal D}_X$-modules  ${\mathcal D}_X
\otimes_{{\mathcal V}^D_0({\mathcal D}_X)}{\mathcal O}_X(D)$ and
${\mathcal D}_X \otimes_{{\mathcal V}^D_0({\mathcal D}_X)}
{\mathcal O}_X$ are dual when $D$ is a free divisor of Spencer type, 
and we have $\mbox{DR}({\mathcal D}_X
\otimes_{{\mathcal V}^D_0({\mathcal D}_X)}{\mathcal O}_X(D) )\cong
\Omega^\bullet_X(\mbox{log}\, D)$. A  generalization of this duality has
been obtained in \cite{CU04}.

\section{Towards an algebraic analogue of {\bf LCT($D$)}\,?}

\label{anal}

In this part, we explain how useful the ${\cal D}_X$-modules are in 
the general study of the condition {\bf LCT($D$)}.

\subsection{Preliminaries}

What can one say about the free divisor $D$ when the morphism
 $\varphi_D$ is an isomorphism? First, we have locally
${\mathcal V}^D_0({\cal D})={\cal O}[\delta_1,\ldots,\delta_n]$ 
where $\{\delta_1,\ldots,\delta_n\}$ is a basis of 
$\mbox{Der}(-\mbox{log\,}D)$ \cite{3}, and ${\cal O}(D)\cong
{\mathcal V}_0^D({\cal D})/
{\mathcal V}_0^D({\cal D})(\delta_1+a_1,\ldots, \delta_n+a_n)$ with 
$\delta_i(h_D)=a_i h_D$, $1\leq i\leq n$. Hence the
morphism $\varphi_D$ is given locally by
\begin{eqnarray*}
  \varphi_D: {\mathcal D}/\tilde{I}^{log} &\longrightarrow &
   {\mathcal O}[1/h_D] \\
   P + \tilde{I}^{log} & \longmapsto & P\cdot\frac{1}{h_D}
\end{eqnarray*}
where $\tilde{I}^{log}\subset {\mathcal D}$ is the left ideal
generated by $\mbox{Ann}_{\cal D}\,1/h_D\cap F_1{\mathcal D}$.
In particular,  it is bijective 
if and only if the two conditions
\begin{description}
   \item[\ \ A($1/h_D$)\,:] the left ideal
  $\mbox{Ann}_{\cal D}\,1/h_D$ of operators
  annihilating $1/h_D$ is generated by operators of order $1$,
   \item[\ \ B($h_D$)\,:] the ${\mathcal D}$-module
  ${\mathcal O}[1/h_D]$ is generated by $1/h_D$,
\end{description}
are satisfied. From a well known result of M. Kashiwara 
(see \cite{11}, Proposition 6.2), this last condition means that 
$-1$ is the only integral root of the
Bernstein polynomial of $h_D$. We recall that the
{\em Bernstein polynomial $b_h(s)$ of $h\in{\mathcal O}$ at the
origin} is the (nonzero) unitary polynomial $b(s)\in{\bf C}[s]$ of
smallest degree which satisfies a functional identity of the form:
$$b(s)h^s=P(s)\cdot h^{s+1}$$
with $P(s)\in{\cal D}[s]={\cal D}\otimes {\mathbf C}[s]$ (see \cite{11}). 
This is an analytic invariant of the ideal $h{\cal O}$.
When $h\in{\mathcal O}$ is not a unit, it is easy to check that $-1$ is
also a root of $b_h(s)$. For example, if $h=x_1^2+\cdots+x_n^2$
then $b_h(s)=(s+1)(s+n/2)$ with the functional identity:
$$(s+1)(s+\frac{n}{2})h^s=\frac{1}{4}\left[\frac{\partial}{\partial x_1}^2 +
\cdots + \frac{\partial }{\partial x_n}^2\,\right]\cdot h^{s+1}\ . $$

\subsection{The conditions {\bf B($h_D$)} and {\bf LCT($D$)}}

\label{BLCT}

 The differential viewpoint above is relevant since it was not
at all clear that for a free divisor, {\bf LCT($D$)} needs the
condition {\bf B($h_D$)}; in particular, every locally weighted 
homogeneous free divisor $D$ satisfies the condition {\bf B($h_D$)} 
at any point (Theorem \ref{CMN} with Theorem \ref{Dual} or 
Proposition \ref{BetLCT}, or more directly \cite{7}, Theorem 5.2). 
Moreover, from the inclusions 
\begin{equation} \label{inclu}
\Omega^\bullet_X(\mbox{log\,$D$})\,\subset\, \Omega^\bullet_X
\otimes ({\cal D}\cdot 1/h_D)\,\subset\, \Omega^\bullet_X\otimes
{\cal O}[1/h_D]  =\Omega^\bullet_X(\star D)
\end{equation}
it appears natural that conditions  {\bf LCT($D$)} and  {\bf B($h_D$)}
are linked for any divisor $D$. As an illustration, we have the
following result:

\begin{prop} \label{BetLCT}
   Let $D\subset X$ be a divisor which satisfies the condition
{\bf LCT($D$)}. Assume that one of the following conditions
is satisfied:
\begin{enumerate}
  \item The divisor $D$ is free except at isolated points.
  \item The divisor $D$ is locally weighted homogeneous.
\end{enumerate}
  Then {\bf B($h_D$)} is satisfied at any point of $D$.
\end{prop}

\begin{demo}
    Firstly we prove the assertion when the condition $1$ is satisfied.
Let $U\subset X$ be a neighborhood of a point 
$m\in D$ such that $D\cap U$ is free at any point different
from $m$. Let $h_D$ be a defining equation of $D$ on $U$.
From Theorem \ref{Dual}, the condition  {\bf B($h_D$)}
is satisfied at any point in $D\cap U-\{m\}$ - since  
$D$ is free at such a point. In particular,
the ${\cal D}_U$-module $\cal C$ in the short exact
sequence
$$0\rightarrow {\cal D}_U\cdot\frac{1}{h_D}
\stackrel{j}{\longrightarrow}
{\cal O}_U(\star D)\longrightarrow {\cal C} \rightarrow 0 $$
is supported at $m$. We just have to prove that $\cal C$ is
zero.  The associated long exact sequence of de Rham cohomology
provides the short exact sequence of ${\bf C}$-vector spaces

$$0 \rightarrow H^n_{DR}({\cal D}_U\cdot\frac{1}{h_D})
\longrightarrow H^n_{DR}({\cal O}_U(\star D))
\longrightarrow H^n_{DR}({\cal C}) \rightarrow 0. $$
In particular, $H^n(j)$ is injective. On the other hand, 
the condition {\bf LCT($D$)} is satisfied; hence,  we deduce from the
inclusions (\ref{inclu}) of de Rham complexes that the 
morphisms $H^i(j):H^i_{DR}({\cal D}_U(1/h_D))\longrightarrow 
H^i_{DR}({\cal O}_U(\star D))$, $0\leq i\leq n$,  are surjective.
Therefore $H^n(j)$ is an isomorphism, that is,  
$H^n_{DR}({\cal C})=0$. From classic results about $\cal D$-modules
supported at a point (see \cite{Mal} for example), $\cal C$ is
necessarily zero and the condition {\bf B($h_D$)} is satisfied at any 
point of $D$.

\stl
  Now, we assume that $D$ is locally weighted homogeneous.
Let us prove the assertion by induction on dimension. 
If $n=2$, then $D$ is locally defined by a (reduced) weighted
homogeneous polynomial in two variables; thus we can
conclude as in the proof of Proposition \ref{ALCTWHIS} 
below. Let us assume that $n\geq 3$ and let $m$ 
denote a point in $D$. From \cite{CMN}, Proposition 2.4,
there exists a neighborhood $U$ of $m$ such that, for
each point $w\in U\cap D$, $w\not=m$, the germ of pair
$(X,D,w)$ is isomorphic to a product 
$({\bf C}^{n-1}\times{\bf C},D'\times {\bf C},(0,0))$ 
where $D'$ is a locally weighted homogeneous divisor of 
dimension $n-2$. Moreover, the condition {\bf LCT($D$)} implies 
that {\bf LCT($D'$)} is satisfied on a neighborhood of the origin 
 (see \cite{CMN}, Lemma 2.2). Let $h_{D'}\in {\cal O}_{{\bf C}^{n-1},0}$
be a local equation of $D'$. By using the induction
hypothesis, the ${\cal D}_{{\bf C}^{n-1},0}$-module 
${\cal O}_{{\bf C}^{n-1},0} [1/h_{D'}]$ is generated by
$1/h_{D'}$, thus so is ${\cal O}_{{\bf C}^{n-1}\times
{\bf C},0}[1/h_{D'}]$. In particular, the condition
{\bf B($h_D$)} is satisfied at any point in $D\cap U-\{m\}$.
We conclude with the first part of the proof.
\end{demo}

In particular,  {\bf LCT($D$)} implies the condition  {\bf B($h_D$)} 
for any divisor $D$ with isolated singularities. For a general
divisor, this question is open.

\subsection{The conditions {\bf A($1/h_D$)} and {\bf LCT($D$)}}

\label{ALCT}

From \cite{T3}, Proposition 1.3, the condition {\bf A($1/h$)} implies
{\bf B($h$)} for any nonzero germ $h\in{\mathcal O}$; 
in particular, {\bf A($1/h_D$)} is a local
analogue of {\bf LCT($D$)} for any Koszul-free divisor $D$, and more
generally, for any free divisor of Spencer type (from Theorem \ref{Dual}; 
see also \cite{CUexp}).
  Is {\bf LCT($D$)} locally equivalent to {\bf A($1/h_D$)} in general?
This question is motivated by the following significant results.
Firstly, this is true for  weighted homogeneous hypersurfaces with an
isolated singularity.

\begin{prop} \label{ALCTWHIS}
    Let $h\in{\bf C}[x]$ be a weighted homogeneous polynomial.
Assume that $h$ defines an isolated singularity at the origin.
Let $D\subset {\bf C}^n$ be the hypersurface defined by $h$.
Then the logarithmic comparison theorem holds for $D$ if
and only if the condition  {\bf A($1/h$)} is satisfied.
\end{prop}

\begin{demo}
  Under our assumptions, the condition  {\bf A($1/h$)} is
in fact equivalent to  {\bf B($h$)} by Theorem \ref{carAIS}.
On the other hand, the polynomial $b_h(s)$ is given
by the formula 
$b_h(s)=(s+1)\prod_{q\in\Pi}(s+|\alpha|+q)$
where $\alpha\in({\bf Q}^{*+})^n$ is the system
of weights such that the degree of $h$ is equal to $1$,
the expression $|\alpha|$  denotes the sum 
$\sum_{i=1}^n\alpha_i\in{\bf Q}^{*+}$, and $\Pi\subset {\bf Q}^+$ is the
set of the degrees of the weighted homogeneous elements
in $A_h={\bf C}[x]/(h'_{x_1},\ldots, h'_{x_n})$ (see \cite{23}, \S 11). 
We recall that $n-2|\alpha|$ is the maximal 
element of $\Pi$;
in particular,  {\bf A($1/h$)} is satisfied if $n=2$ and so is  {\bf LCT($D$)}
(Theorem \ref{CCMN}). Moreover, the set $\Pi$ is symmetric
about $(n/2)-|\alpha|$; 
hence, we deduce easily that {\bf B($h$)} is
equivalent to the last condition of Theorem \ref{hachem} when
$n\geq 3$. This completes the proof.
\end{demo}

Moreover, the condition {\bf B($h_D$)} is satisfied by any
hyperplane arrangement (A. Leykin \cite{Walt}, Theorem 5.1), 
and the condition 
{\bf A($1/h_D$)} is true for the union of a generic arrangement 
with a hyperbolic arrangement \cite{T3}; this agrees with Terao's 
conjecture (see \S \ref{arrang}). 
  The general problem is still open,
and  condition {\bf A($1/h$)}  may only be  necessary. 
A difficulty is the lack of families of divisors which
satisfy the condition {\bf LCT($D$)}.

\subsection{The condition {\bf A($1/h$)}}

\label{AsurH}

Let $h\in{\mathcal O}$ be a nonzero germ such that $h(0)=0$.
We give here some results about the meaning of the condition
{\bf A($1/h$)} (see \cite{T3}).

\stl

 First, we have the following easy criterion:
\begin{lemm} \label{crit}
  Let $h\in{\mathcal O}$ be a nonzero germ such that $h(0)=0$.
Assume that the following conditions are satisfied:
\begin{description}
   \item[\ \ H($h$)\,:] $h$ belongs to the ideal of its partial derivatives;
   \item[\ \ B($h$)\,:] $-1$ is the smallest integral root of $b_h(s)$;
   \item[\ \ A($h$)\,:] the ideal $\mbox{\em Ann}_{\cal D}\,h^s$ is
  generated by operators of order $1$.
\end{description}
Then the ideal $\mbox{\em Ann}_{\cal D}\,1/h$ is generated by
operators of order $1$.
\end{lemm}

\begin{demo}
By Euclidean division, we have also a decomposition
$$\mbox{Ann}_{{\cal D}[s]}\,h^s={\mathcal D}[s](s-v)+
{\mathcal D}[s]\mbox{Ann}_{\cal D}\,h^s$$
where $v$ is a vector field such that $v(h)=h$ (condition {\bf H($h$)}).
Moreover, under the condition {\bf B($h$)}, the ideal
$\mbox{Ann}_{\cal D}\,1/h$ is obtained by fixing $s=-1$ in
a system of generators of $\mbox{Ann}_{{\cal D}[s]}\,h^s$.
\end{demo}

Reciprocally, what does remain true? We recall that
 the condition {\bf A($1/h$)} always implies {\bf B($h$)}. 
On the other hand, does {\bf A($1/h$)} imply {\bf H($h$)}? 
This is true for isolated singularities \cite{T1}, 
Koszul-free germs, and suspensions of unreduced 
plane curve $z^N+g(x_1,x_2)$ (see \cite{T3}); this
question is still open. 
Finally, the condition {\bf A($1/h$)} does not imply
{\bf A($h$)} in general. Indeed, Calder\'on's example
$h=x_1x_2(x_1+x_2)(x_1+x_2x_3)$ satisfies {\bf LCT($D$)},
{\bf A($1/h$)}, {\bf B($h$)}, {\bf H($h$)}
and not {\bf A($h$)} (see \cite{3}, \cite{CCMN}, 
\cite{CN}, \cite{CUE}, \cite{T3}). 
Meanwhile, condition {\bf A($h$)} is
not unrealistic, since we have the following
characterization of {\bf A($1/h$)} for Koszul-free germs:

\begin{theo} {\em \cite{T3}}    \label{carkosfree}
   Let $h\in{\mathcal O}$ be a Koszul-free germ.
Then the left ideal $\mbox{{\em Ann}}_{\mathcal D}\,1/h$ is
generated by operators of order one  if and only if
the conditions {\bf H($h$)}, {\bf B($h$)} and {\bf A($h$)}
are satisfied.
\end{theo}

Moreover,  condition {\bf A($h$)} is satisfied 
when $h$ defines an isolated singularity (see below).
Thus, we have 

\begin{theo} {\em \cite{T1}} \label{carAIS}
   Let $h\in{\mathcal O}$ be a germ of holomorphic
function defining an isolated singularity.
Then the ideal $\mbox{{\em Ann}}_{\mathcal D}\,1/h$ is
generated by operators of order one  if and only if
the germ $h$ is weighted homogeneous  and the condition
{\bf B($h$)} is satisfied.
\end{theo}

In fact, the condition  {\bf A($h$)} may be considered almost
as a geometric condition. Indeed, the following condition implies
 {\bf A($h$)}:
\begin{description}
   \item[\ \ W($h$)\,:] the relative conormal space
$W_h$ is defined by linear equations in $\xi$
\end{description}
since $W_h=\overline{\{(x,\lambda\, \mbox{d}h)\,:\,
\lambda\in{\mathbf C}) \}}\subset T^*{\mathbf C}^n$ is the
characteristic variety of ${\mathcal D}h^s$ \cite{11}. For
example, {\bf W($h$)} is true for hypersurfaces with an isolated
singularity \cite{23} and for locally weighted homogeneous 
free divisors \cite{CN}. 

\bigskip

\stl

\noindent{\bf\Large Appendix}

\stl

The following diagram summarizes the
expected relations between the conditions
studied in this note:

\begin{figure}[ht]
\begin{center}
\psfrag{P1}{\S \ref{free}}
\psfrag{P2}{\S \ref{AsurH}}
\psfrag{P3}{\S  \ref{ALCT}}
\psfrag{P4}{\S \ref{arrang}}
\psfrag{P5}{\S \ref{BLCT}}
\psfrag{P6}{\S \ref{ALCT}}
\psfrag{R1}{\cite{Walt}}
\psfrag{R2}{\cite{T3} }
\psfrag{LCT}{{\bf LCT($D$)}}
\psfrag{A}{{\bf A($1/h_D$)}}
\psfrag{B}{{\bf B($h_D$)}}
\psfrag{H}{{\bf H($h_D$)}}
\includegraphics[height=6cm]{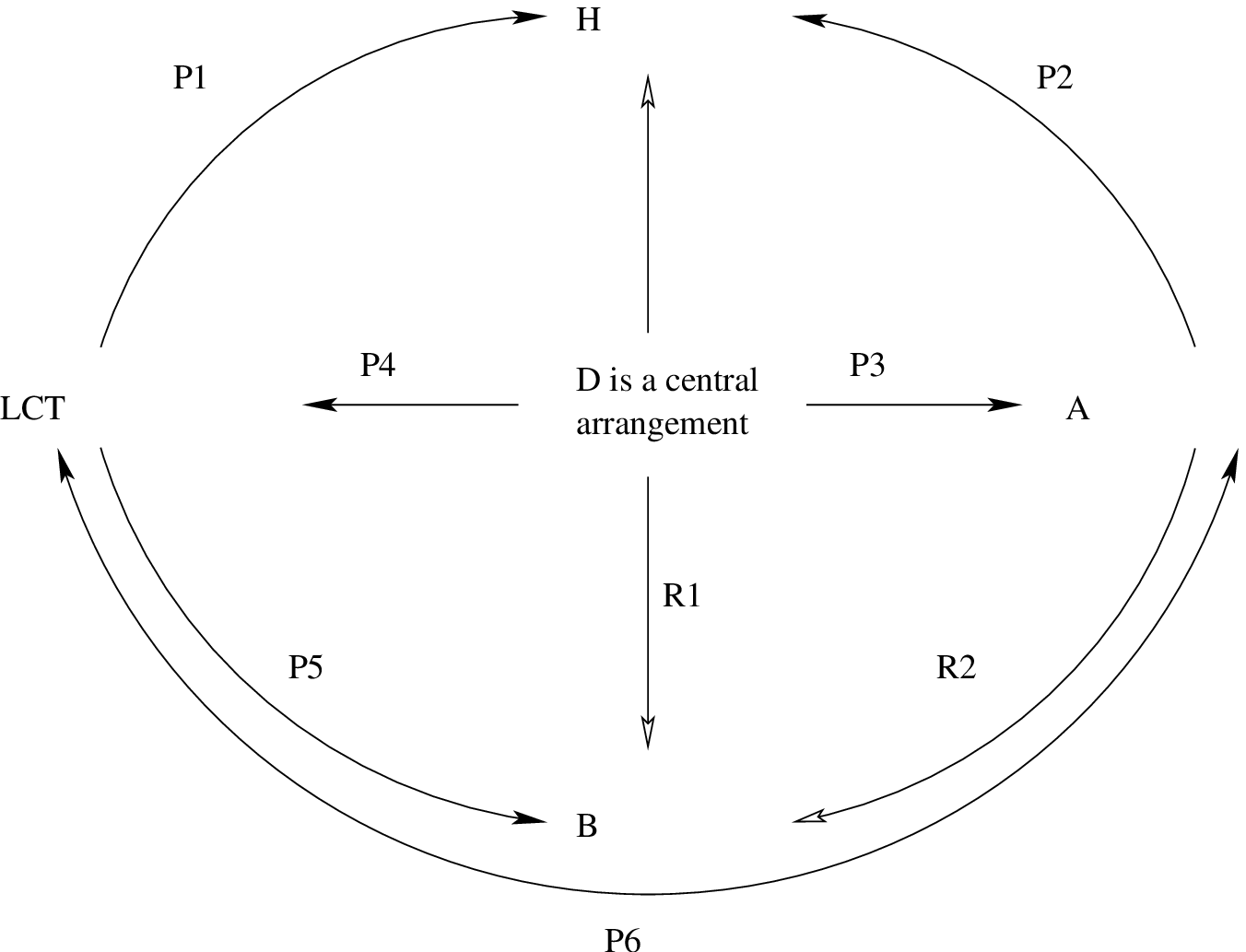}
 
\end{center}
\end{figure}

Here $D$ denotes a general divisor, and
the solid headed arrows represent the open
questions.

\end{document}